\documentclass[10pt,a4paper]{article}
\usepackage[latin1]{inputenc}
\usepackage{fullpage}
\usepackage{mathrsfs}
\usepackage{amsmath}
\usepackage{amsfonts}
\usepackage{amssymb}
\usepackage{amsthm}
\usepackage{graphicx}
\usepackage{epstopdf}
\usepackage[affil-it]{authblk}
\graphicspath{C:/Users/evanc/OneDrive/Documents/The distributional differintegral}
\DeclareGraphicsExtensions{.eps,.pdf,.png,.jpg}
\author{Evan Camrud%
	\thanks{Electronic address: \texttt{ecamrud@iastate.edu}; Corresponding author. Current address: Mathematics Department, Iowa State University, Ames, IA, 50014} }
\affil{Mathematics Department, Concordia College, Moorhead, MN, 56562}
\title{A Novel Approach to Fractional Calculus: Utilizing Fractional Integrals and Derivatives of the Dirac Delta Function}
\date{ }

\theoremstyle{definition}
\newtheorem{Definition}{Definition}

\begin{document}
	
	\maketitle
	

	\textbf{Abstract.}\footnote{Fractional calculus, fractional differential equations, integral transforms, operations with distributions, special functions. \textbf{2010 Mathematics Subject Classification}.  Primary 26A33; Secondary 34A08, 46F10, 33E99} {While the definition of a fractional integral may be codified by Riemann and Liouville, an agreed-upon fractional derivative has eluded discovery for many years. This is likely a result of integral definitions including numerous constants of integration in their results. An elimination of constants of integration opens the door to an operator that reconciles all known fractional derivatives and shows surprising results in areas unobserved before, including the appearance of the Riemann Zeta function and fractional Laplace and Fourier Transforms. A new class of functions, known as Zero Functions and closely related to the Dirac delta function, are necessary for one to perform elementary operations of functions without using constants. The operator also allows for a generalization of the Volterra integral equation, and provides a method of solving for Riemann's ``complimentary'' function introduced during his research on fractional derivatives.}

	\section{Introduction}
	
	The concept of derivatives of non-integer order, commonly known as fractional derivatives, first appeared in a letter between L'Hopital and Leibniz in which the question of a half-order derivative was posed \cite{leibniz1962letter}. In recent years, the research has found footholds in many areas of study, including applications in polymers, quantum mechanics, group theory, wave theory spectroscopy, continuum mechanics, field theory, biophysics, statistics, and Lie theory \cite{bapna2012application,herrmann2014fractional,hilfer2000applications,hilfer2008threefold,machado2011and}. Many formulations of fractional derivatives have appeared over the centuries, such as the Riemann-Liouville, Caputo, Hadamard, Erdelyi-Kober, Grunwald-Letnikov, Marchaud, and Riesz, but one would expect an ``ultimate'' definition to emerge out of the many \cite{kilbas2006theory,miller1993introduction,oldham1974fractional,podlubny1998fractional,samko1993fractional}. This ``ultimate'' has seemingly eluded discovery, and one is forced to choose a so-called ``best derivative for the job'', depending on how a particular definition relates to the research at hand. 
	
	That is not to say, however, that the perfect definition does not exist. It seems likely that one could expect the following to be true:
	\begin{enumerate}
		\item $\frac{d^\alpha}{dx^\alpha}x^n=\frac{\Gamma(n+1)}{\Gamma(n+1-\alpha)}x^{n-\alpha}$ for $n\geq0$ and $\alpha\leq n+1$,
		\item $\frac{d^\alpha}{dx^\alpha}e^{\lambda x}=\lambda^\alpha e^{\lambda x}$,
		\text{which, assuming the derivative is linear, implies}
		\item $\frac{d^\alpha}{dx^\alpha}\sin(\lambda x)=|\lambda|^\alpha \sin(\lambda x+\frac{\pi}{2}\alpha)$,
		\text{and}
		\item $\frac{d^\alpha}{dx^\alpha}\cos(\lambda x)=|\lambda|^\alpha \cos(\lambda x+\frac{\pi}{2}\alpha)$.
	\end{enumerate}
	This is a result of noticing the patterns of traditional derivatives, and interpolating their properties. Thus far, no proposed definition satisfies all four of the above in all cases, and indeed there is much debate as to whether the above are truly the ``correct'' interpolations of their respective patterns.
	
	The discrepancies inherent in fractional derivative definitions are likely due to the fact that nearly all fractional derivatives are instead based on generalizing repeated \textit{integration}. This brings up many questions such as, ``What should the upper-and-lower limits of integration be?'', or ``Should there be terms added to the end to cancel out abnormalities?''
	
	Even so, a commonly used definition for the fractional derivative is the Riemann-Liouville definition, which is a generalization of Cauchy's formula for repeated integration: 
	\begin{equation}\frac{1}{\Gamma (\alpha)} \int_{c}^{x} f(\tau)(x-\tau)^{\alpha-1} d\tau,\end{equation}
	with $c$ as an arbitrary integration limit. This, however, is by its nature a \textit{fractional integral}. To make the fractional integral into a derivative, a full derivative of the fractional integral is taken. This definition introduces surprising results, such as the fractional derivative of a constant not being constant. Caputo eliminated this ``abnormality'' by adding a small term onto the end which would subtract whatever a constant evaluated in the Riemann-Liouville definition, leaving zero.
	
	Another very popular definition, the Grunwald-Letnikov fractional derivative arises from a binomial generalization of repeated limit-based derivatives,
	\begin{equation}\frac{d^\alpha}{dx^\alpha}f(x)=\lim_{h\to0}\sum_{m=0}^\infty \frac{(-1)^{m+\alpha}}{h^\alpha}\binom{\alpha}{m}f(x+mh).\end{equation}
	This derivative is also special in that it can provide results for complex values of $\alpha$. Acting in this manner upon the exponential function allows for a wide range of use within harmonic analysis, wavelet theory, and other branches of mathematics that deal with Fourier series \cite{barbosa2007analysis,wu2009wavelet}.
	
	While they do not satisfy all of the aforementioned conjectured results,  the Riemann-Liouville and Grunwald-Letnikov derivatives indeed satisfy the four properties in the \textit{following} definition, which may be taken as the definition of a fractional derivative, as defined by Ortigueira and Machado \cite{ortigueira2015fractional}:
	
	\begin{Definition}
		
		Let $\alpha \in [0,1]$.  An operator $D^\alpha$ is a fractional differential operator if it satisfies the following four properties:
		
		\begin{enumerate}
			\item $\text{Linearity: } D^\alpha(af + bg) = a D^\alpha(f) + b  D^\alpha(g) \text{ for all }a, b \in \mathbb{C} \text{ and } f, {g} \in \text{Dom}(D^\alpha), \text{ where } \text{Dom}(D^\alpha)$
			
			$\text{ is the domain of the operator } D^\alpha$
			\item $D^0[f]=f \text{ for all functions } f$
			\item $D^1[f]=f'\text{ for all } f \in \text{Dom}(D^1)$
			\item $\text{The Index Law: } D^\beta D^\alpha [f]=D^{\beta+\alpha}[f] \text{ for all } f \in \text{Dom}(D^\beta \circ D^\alpha) \cap \text{Dom}(D^{\beta+\alpha})$.
		\end{enumerate}
		
	\end{Definition}
	
	Satisfying the above definition is good, but not quite good enough to ``win-out'' against all other forms of fractional derivatives. And just so, newer definitions are arising that bend these rules so that other rules may be met instead \cite{anderson2015newly}. This rule-bending allowed for the first non-linear conformable fractional derivative to be proposed just two years ago \cite{camrud2017conformable}.
	
	The manuscript which follows eliminates constant functions, and in doing so changes the nature of the spaces which the fractional derivative behaves, forcing its domain into generalized function spaces. This in turn allows the definition to give results in terms of \textit{distributional} derivatives and integrals, and even changes the notion of integrals themselves.
	
	In comparing the results to past distributional fractional derivatives, these results may  be a reconciliation of what was proposed in \cite{stojanovic2011generalized}, namely that the Riemann-Liouville derivative operating in a distributional sense does not produce an integer-valued distributional derivative. The results obtained in \cite{stojanovic2011generalized}, when paired with the notion of Zero Functions proposed here, give precisely the integer-valued distributional derivatives.
	
	This definition allows derivation and integration of complex powers, and does so with a single definition between both derivation and integration, making the operator to the negative power the inverse. The definition also allows for the construction of fractional integral transforms, the solving of fractional differential equations with an arbitrary number of initial/boundary conditions. What follows is the introduction of the distributional differintegral, and an overview of its many properties.
	
	\section{Elimination of Nonzero Constants from Allowed Functions of Differintegration}
	
	\textbf{Remark.} This paper uses the terms ``antiderivative'' and ``integral'' interchangeably.
	
	Derivatives and antiderivatives are \textit{not} inverses of one another. Considering the function $f(z)=1$ and the integral and differential operators $J$ and $D$ respectively,
	
	\begin{equation}D^n J^n [1]=1 \text{ but } J^n D^n[1]=\sum_{k=0}^{n-1} c_k x^k.\end{equation}
	
	This is why $n$ boundary/initial conditions are necessary for a differential equation of order $n$. But this is also why different integral definitions of fractional derivatives (or even different bounds of integration on the \textit{same} definition) yield drastically different results.
	
	As recently shown by \cite{labora2018index}, the secret to fractional differentiation lies in eliminating nonzero constants from ``allowed'' functions of differentiation/integration. This is a result of the differential operator losing its bijectivity on a domain containing these functions. Thus, it cannot be invertible. Instead, it is apt to allow the \textit{distributions} of the form
	
	\begin{equation}f(z)=Cz^0=C\big[H(z)+H(-z)\big]\end{equation}
	
	(where $H(z)$ is the Heaviside step function) which are equal to constant functions \textit{almost everywhere}, but remain undefined at $z=0$.
	
	One must also use the following identity given by the power rule for derivatives (but not the limit definition):
	
	\begin{equation}\frac{d}{dz}z=z^0=H(z)+H(-z).\end{equation}
	
	To avoid confusion with past definitions of antiderivatives, however, the new system applied uses the terminology ``inverse derivative'', along with the operator $\frac{d^{-1}}{dz^{-1}}$ for the inverse derivative of a function of independent variable $z$.
	
	Since nonzero constant functions are no longer allowed, it may be enforced that
	
	\begin{equation}\frac{d^{-1}}{dz^{-1}} 0=0\text{ and }\frac{d^{-1}}{dz^{-1}}z^0=z.\end{equation}
	Note that one may relate an inverse derivative to an integral by
	
	\begin{equation}\int f(z)dz=\frac{d^{-1}}{dz^{-1}}f(z)+C.\end{equation}
	
	\begin{Definition}
		Let $X(\Omega)$ be a linear generalized function space containing the monomials. Define $X^{-c}(\Omega)$ to be the generalized function space containing precisely the image of the linear mapping, $T:X(\Omega)\to X^{-c}(\Omega)$ and all derivatives (in the sense of $X(\Omega)$) of the image of $T$ on $X(\Omega)$.
		
		The mapping $T$ is defined as follows
		
		\begin{equation}T\big[f(z)\big]=\begin{cases}
		f(z)&\text{ if }f\text{ is not a constant function}\\
		
		f(z)\cdot z^0& \text{ if }f\text{ is a constant function}
		\end{cases}\end{equation}
		for all $f\in X(\Omega)$, and where $z^0$ is assumed to have the properties explained above in this section.
		
		The space $X^{-c}(\Omega)$ is called a generalized function space \underline{with trivial constant}.
	\end{Definition}
	
	From this point onwards, only derivatives and inverse derivatives in $X^{-c}(\Omega)$ spaces are considered.
	
	\subsection{Definition of the Zero Function}
	
	Constant functions have been eliminated and replaced with functions equal to constants almost everywhere. To utilize these, it is reiterated that
	
	\begin{equation}\frac{d}{dz}z=z^0=H(-z)+H(z).\end{equation}
	This results in very novel cases, such that
	
	\begin{equation}\frac{d}{dz}z^0=0\cdot z^{-1}\simeq\delta(z)-\delta(-z),\end{equation}
	where $\delta(z)$ is the Dirac delta function. This notation is used to give an intuitive understanding of the shape of this derivative as an ``odd'' function. The distributional derivative of $z^0=H(-z)+H(z)$ is indeed $\delta(z)-\delta(-z)$ however the \textit{approximately-equal} is used to indicate the use of the \textit{inverse distributional derivative}, which is slightly different than a distributional integral.
	
	Herein is proposed the Zero Function:
	
	\begin{Definition}
		The Zero Function $\emptyset(z)$ is defined as
		
		\begin{equation}\emptyset(z)=\frac{d}{dz}z^0\simeq\delta(z)-\delta(-z).\end{equation}
		
	\end{Definition} 
	Utilizing inverse derivatives one obtains 
	
	\begin{equation}\frac{d^{-1}}{dz^{-1}}\emptyset(z)=z^0.\end{equation}
	
	One may continue to take derivatives of the Zero Function:
	
	\begin{equation}\emptyset^{(n-1)}(z)=\frac{d^n}{dz^n}z^0=H^{(n)}(z)+(-1)^n H^{(n)}(-z)\simeq\delta^{(n-1)}(z)-\delta^{(n-1)}(-z).\end{equation}
	
	It results that one must reconcile the idea of the Zero Function to an operator that is \textit{not} an integral, but rather an inverse derivative.
	
	\section{Interpolating the differintegral of the Heaviside step function}
	
	In the sense of distributions, it is well known that $H'(x)=\delta(x)$, and also that $\int_{-\infty}^xH(t)dt=xH(x)$.
	
	Generalizing these with differential and integral operators one obtains
	
	\begin{equation}D^n[H(x)]=\delta^{(n-1)}(x),\text{ and }J^n[H(x)]=\frac{x^n}{n!}H(x).\end{equation}
	Interpolating these results one may argue that
	
	\begin{equation}J^\alpha[H(x)]=\frac{x^{\alpha}}{\Gamma(1+\alpha)}H(x).\end{equation}
	
	Observe now the integral to a negative integer power
	
	\begin{equation}J^{-n}[H(x)]=\frac{x^{-n}}{\Gamma(1-n)}H(x).\end{equation}
	To see how this function behaves as a distribution, one must act it upon a test function, $\phi$.
	
	\begin{equation}
	\begin{split}
	&\int_{\mathbb{R}} \frac{(x-y)^{-n}}{\Gamma(1-n)}H(x-y) \phi(y)dy=\frac{1}{\Gamma(1-n)}\int_{\mathbb{R}}\phi(y)H(x-y)(x-y)^{-n} dy=\frac{-n}{\Gamma(1-n)}\int_{\mathbb{R}}\phi'(y)H(x-y)(x-y)^{1-n} dy\\
	&=\frac{-n(1-n)}{\Gamma(1-n)}\int_{\mathbb{R}}\phi''(y)H(x-y)(x-y)^{2-n} dy=...\simeq\frac{\Gamma(1-n)}{\Gamma(1-n)}\int_{\mathbb{R}}\phi^{(n+1)}(y)H(x-y) dy=\phi^{(n)}(x)
	\end{split}
	\end{equation}
	where the equality after the ellipsis is considered in the sense of residues.
	
	At each equality above, integration by parts was used, but because test functions are only nonzero on a compact set, they vanish at infinity. Indeed, this realization that the derivative of the Heaviside step function may be fractionalized to positive and negative powers	allows one to construct the distributional differintegral.
	
	\section{Definition of the distributional differintegral}
	
	To utilize any notion of a ``distributional'' differintegral, one must understand what spaces allow a distributional derivative to be taken. Thus, a function space must contain a dense subset of functions analogous to test functions to make any progress in defining a distributional differintegral. Herein are defined the spaces to consider:
	
	\begin{Definition}
		Let $X(\Omega)$ be a normed space of functions defined on $\Omega$ which includes a definition for a derivative (and as such an inverse derivative). Then define $X(\Omega)$ to be a \underline{distributional function space} if there exists a dense (with respect to $\|\cdot\|_{X(\Omega)}$) subset $X_t(\Omega)\subseteq X(\Omega)$ such that for all $f\in X_t(\Omega)$,
		
		\begin{enumerate}
			\item $\frac{d^n}{dz^n}f(z)=f^{(n)}(z)$ exists for all $n\in\mathbb{N}$ and $z\in \Omega$.
			\item There exists a compact subset $K\subseteq\Omega$, such that $f^{(n)}(z)=0$ for all $z\in \Omega\backslash K$ and for all $n\in\mathbb{N}$.
		\end{enumerate}
		
		Define each element $f\in X_t(\Omega)$ to be a \underline{test function} of the space $X(\Omega)$.
		
	\end{Definition}
	One may recognize the above as a generalization of the test functions in many spaces of real-valued functions.
	
	Just as derivatives and integrals have often been defined, an operator is used to represent the distributional differintegral. This operator is an elongated section sign, \S, with the variable of differintegration subscripted, and the power of differintegration superscripted. The motivation for this is that the symbol is close to that of both $\int$ and $\delta$.
	
	The definition of the distributional differintegral is thus presented.
	
	\begin{Definition}
		
		Let $f\in X_t(\Omega)\subseteq X^{-c}(\Omega)$ a distributional function space with trivial constant, and let $z_0\in\partial \Omega$ (note $z_0$ may be infinite).
		
		For $\alpha\in\mathbb{C}$, the $\alpha^{th}$ distributional differintegral of $f(z)$, with respect to the variable $z$, is
		
		\begin{equation}
		\begin{split}
		\text{\LARGE\S}_z ^\alpha f(z)&=\frac{1}{\Gamma(\alpha)}\bigg(\frac{d^{-1}}{d\zeta^{-1}} f(\zeta)(z-\zeta)^{\alpha-1}\Big\vert_{\zeta=z}\bigg)\\
		&=\frac{1}{\Gamma(\alpha)}\int_{z_0}^z f(\zeta)(z-\zeta)^{\alpha-1}d\zeta\\	
		&=\frac{1}{\Gamma(\alpha)}\oint_{\gamma} f(\zeta)(z-\zeta)^{\alpha-1}H(z-\zeta)d\zeta
		\end{split}
		\end{equation}
		where $\gamma$ is a simple closed curve in $\Omega$ containing the points $z_0$ and $z$, with $H(z-\zeta)$ becoming the real-valued Heaviside step function when $\gamma$ is parameterized by a real variable.\footnote{That is, if $\gamma(t):[t_0,t_1]\to \Omega$ with $\gamma(t_0)=\gamma(t_1)=z_0$ and $\gamma(t')=z$ for $t'\in[t_0,t_1]$, then $H(z-\zeta)=H\big(z-\gamma(t)\big)=1$ if $t\in[t_0,t')$ and $H(z-\zeta)=H\big(z-\gamma(t)\big)=0$ if $t\in(t',t_1]$.}
	\end{Definition}
	
	\textbf{Remark.} When working with functions of a single real variable, the simple closed curve $\gamma$ may be thought of as the real-part of a circle on the complex Riemann-sphere. This includes the most common case where the curve becomes the real line $[-\infty,\infty]$.
	
	When the integrals in the definition above do not converge, one may evaluate the integral for values of $\alpha$ which converge and perform \textit{analytic continuation} (with respect to $g_z(\alpha)=\text{\Large\S}_z ^\alpha f(z)$ as a function of $\alpha$) to give valid results for all $\alpha\in\mathbb{C}$. Since $\frac{1}{\Gamma(\alpha)}$ is entire, this function will be well-defined.
	
	Any of the three definitions above are equivalent. It is sometimes helpful to understand the first (inverse derivative) definition, noting that this definition does not compute integrals but inverse derivatives. Upon first inspection, however, one may see the extreme similarity between this definition and that of the Riemann-Liouville fractional integral: 
	
	\begin{equation}_{c}J_x^\alpha[f(x)]=\frac{1}{\Gamma (\alpha)} \int_{c}^{x} f(\tau)(x-\tau)^{\alpha-1} d\tau.\end{equation}
	
	\textbf{Remark.} The definitions are \textit{equivalent} for real-valued test functions, $\alpha>0$, and a lower integration bound of $-\infty$.
	
	An inverse derivative definition (instead of only an integral definition) sets no limits on convergence. This allows extension values of $\alpha$ to all complex numbers as explained above, instead of just the positive real numbers.
	
	It is important to insist that the Zero Function be used, as it defines an inverse derivative without a nonzero constant of integration. Unless a derivative explicitly denotes the value of the constant of integration (by means of an $n^{th}$ derivative of the Zero Function), there \textit{is no constant of integration}.
	
	It is important to note here that $\alpha\in\mathbb{Z}^+$ computes the $\alpha^{th}$ inverse derivative, and $\alpha\in\mathbb{Z}^-$ computes the $\alpha^{th}$ derivative, while $\alpha=0$ is the identity operation.
	
	\subsection{The distributional differintegral for non-test functions and distributions}
	
	Since the above definition holds only for test functions, one must extend the distributional differintegral to other functions, as well as distributions. Since it was required that the test functions be dense in the function space, one may define the differintegral for non-test functions as follows:
	
	\begin{Definition}
		Let $f\in X^{-c}(\Omega)$, a distributional function space with trivial constant. Since $X_t(\Omega)$ is dense in $X^{-c}(\Omega)$ there exists $\{\phi_n\}_{n=1}^\infty$ such that $\phi_n\in X_t(\Omega)$ for all $n$ and $\phi_n\to f$ in $\|\cdot\|_{X^{-c}(\Omega)}$. For $\alpha\in\mathbb{C}$, the $\alpha^{th}$ distributional differintegral of $f(z)$, with respect to the variable $z$, is
		
		\begin{equation}\text{\LARGE\S}_z^\alpha f(z)=\lim_{n\to\infty}\text{\LARGE\S}_z^\alpha \phi_n(z).\end{equation}
	\end{Definition}
	
	In regards to distributions, as with distribution theory, the distributional differintegral is defined as follows.
	
	\begin{Definition}
		Let $T\in X_t(\Omega)^*$, the space of continuous linear functionals on $X_t(\Omega)$. For $\alpha\in\mathbb{C}$, the $\alpha^{th}$ distributional differintegral of $T$, with respect to the variable $z$, is
		
		\begin{equation}\text{\LARGE\S}_z^\alpha T\big[(\circ)\big]=e^{-i\pi\alpha}\cdot T\bigg[\text{\LARGE\S}_z^\alpha (\circ)\bigg]\end{equation}
		where $(\circ)$ is a placeholder for test functions.
	\end{Definition}
	Notice that if $T_f$ is a regular distribution, with $f$ a locally integrable function, and $T\big[(\circ)\big]=\int_{\Omega}\big[f(z)\cdot(\circ)\big]dz$, then $\text{\Large\S}_z^\alpha T_f=T_{\text{\large\S}_z^\alpha f}$.
	
	\section{The distributional differintegral is a fractional derivative}
	
	Arising from Ortigueira and Machado's definition of a fractional differential operator, it should be shown that the distributional differintegral satisfies all four properties. The following are proofs only for test functions of a single real variable, but as they are dense and the definition for other functions is dependent on that of test functions, these are the only proofs given. Proofs for functions in general spaces may be researched in the future.
	
	\textbf{Proof.}
	
	\begin{enumerate}
		\item Linearity: $\text{\LARGE\S}_z ^\alpha \Big[\lambda f(z)+\mu g(z)\Big]=\lambda\text{\LARGE\S}_z ^\alpha f(z)+\mu\text{\LARGE\S}_z ^\alpha g(z)$.
		
		\begin{equation}
		\begin{split}
		\text{\LARGE\S}_z ^\alpha \Big[\lambda f(z)+\mu g(z)\Big]&=\frac{1}{\Gamma(\alpha)}\bigg(\frac{d^{-1}}{d\zeta^{-1}} \big(\lambda f(\zeta)+\mu g(\zeta)\big)(z-\zeta)^{\alpha-1}\Big\vert_{\zeta=z}\bigg)\\
		&=\frac{1}{\Gamma(\alpha)}\bigg(\frac{d^{-1}}{d\zeta^{-1}} \lambda f(\zeta)(z-\zeta)^{\alpha-1}+\mu g(\zeta)(z-\zeta)^{\alpha-1}\Big\vert_{\zeta=z}\bigg)\\
		&=\frac{1}{\Gamma(\alpha)}\Bigg[\bigg(\frac{d^{-1}}{d\zeta^{-1}} \lambda f(\zeta)(z-\zeta)^{\alpha-1}\Big\vert_{\zeta=z}\bigg)+\bigg(\frac{d^{-1}}{d\zeta^{-1}}\mu g(\zeta)(z-\zeta)^{\alpha-1}\Big\vert_{\zeta=z}\bigg)\Bigg]\\
		&=\frac{1}{\Gamma(\alpha)}\Bigg[\lambda\bigg( \frac{d^{-1}}{d\zeta^{-1}} f(\zeta)(z-\zeta)^{\alpha-1}\Big\vert_{\zeta=z}\bigg)+\mu\bigg(\frac{d^{-1}}{d\zeta^{-1}} g(\zeta)(z-\zeta)^{\alpha-1}\Big\vert_{\zeta=z}\bigg)\Bigg]\\
		&=\lambda\frac{1}{\Gamma(\alpha)} \bigg( \frac{d^{-1}}{d\zeta^{-1}} f(\zeta)(z-\zeta)^{\alpha-1}\Big\vert_{\zeta=z}\bigg)+\mu\frac{1}{\Gamma(\alpha)}\bigg(\frac{d^{-1}}{d\zeta^{-1}} g(\zeta)(z-\zeta)^{\alpha-1}\Big\vert_{\zeta=z}\bigg)\\
		&=\lambda\text{\LARGE\S}_z ^\alpha f(z)+\mu\text{\LARGE\S}_z ^\alpha g(z).
		\end{split}
		\end{equation}
		
		\item For all test functions $f(z)$, $\text{\LARGE\S}_z ^0 f(z)=f(z)$.
		
		\begin{equation}\text{\LARGE\S}_z ^0 f(z)=\int_{\mathbb{R}} f(\zeta)\frac{(z-\zeta)^{-1}}{\Gamma(0)}H(z-\zeta)d\zeta=\int_{\mathbb{R}} f(\zeta)\delta(z-\zeta)d\zeta=f(z).\end{equation}
		
		(See section 3 above for further explanation.)
		
		\item For all test functions $f(z)$, $\text{\LARGE\S}_z ^{-1} f(z)=f'(z)$.
		
		\begin{equation}\text{\LARGE\S}_z ^{-1} f(z)=\int_{\mathbb{R}}f(\zeta)\frac{(z-\zeta)^{-2}}{\Gamma(-1)}H(z-\zeta)d\zeta=\int_{\mathbb{R}} f(\zeta)\delta'(z-\zeta)d\zeta=f'(z).\end{equation}
		
		(See section 3 above for further explanation.)
		
		\item The Index Law: For all test functions $f(z)$, $\text{\LARGE\S}_z ^\alpha \bigg[\text{\LARGE\S}_z ^\beta f(z)\bigg]=\text{\LARGE\S}_z ^{\alpha+\beta}f(z)$.\\
		
		To solve this, one must use the beta function,
		
		\begin{equation}\int_0^1 u^{\alpha-1}(1-u)^{\beta-1}du=\frac{\Gamma(\alpha)\Gamma(\beta)}{\Gamma(\alpha+\beta)},\end{equation}
		as well as the Dirichlet formula, given by \cite{whittaker1952analysis}, but in the form necessary for the proof, 
		
		\begin{equation}\int_{-\infty}^z (z-\zeta)^{\alpha-1}d\zeta\int_{-\infty}^{\zeta}f(\phi)(\zeta-\phi)^{\beta-1}d\phi=\int_{-\infty}^zf(\phi)\bigg[\int_{\phi}^{z}(z-\zeta)^{\alpha-1}(\zeta-\phi)^{\beta-1}d\zeta\bigg]d\phi.\end{equation}
		Therefore one obtains
		
		\begin{equation}\text{\LARGE\S}_z^\alpha \bigg[\text{\LARGE\S}_z^\beta f(z)\bigg]=\frac{1}{\Gamma(\alpha)\Gamma(\beta)}\int_{-\infty}^{z}f(\phi)\bigg[\int_\phi^z(z-\zeta)^{\alpha-1}(\zeta-\phi)^{\beta-1}d\zeta \bigg]d\phi.\end{equation}
		
		The inner integral, $k(z,\phi)=\int_\phi^z(z-\zeta)^{\alpha-1}(\zeta-\phi)^{\beta-1}d\zeta$, may be interpreted as the kernel of the external convolution integral. With the substitution $u=\frac{\zeta-\phi}{z-\phi}$, which leads to $\zeta=\phi+u(z-\phi)$ and $d\zeta=(z-\phi)du$, one obtains
		
		\begin{equation}
		\begin{split}
		k(z,\phi)&=\int_{\phi}^{z}(z-\zeta)^{\alpha-1}(\zeta-\phi)^{\beta-1}d\zeta\\	
		&=(z-\phi)^{\alpha+\beta-1}\int_0^1 u^{\alpha-1}(1-u)^{\beta-1}du\\
		&=\frac{\Gamma(\alpha)\Gamma(\beta)}{\Gamma(\alpha+\beta)}(z-\phi)^{\alpha+\beta-1}.
		\end{split}
		\end{equation}
		
		Thus one may complete the proof,
		
		\begin{equation}\text{\LARGE\S}_z ^\alpha \bigg[\text{\LARGE\S}_z ^\beta f(z)\bigg]=\frac{1}{\Gamma(\alpha+\beta)}\int_{-\infty}^z f(\phi)(z-\phi)^{\alpha+\beta-1}d\phi=\text{\LARGE\S}_z ^{\alpha+\beta}f(z).\end{equation}
		
		\begin{flushright}
			$\square$
		\end{flushright}
		
	\end{enumerate}
	
	\section{Specific relationship to the Riemann-Liouville definition}
	
	In Riemann's initial, posthumous publication of fractional calculus, his definition was as follows:
	
	\begin{equation}\frac{d^{-\alpha}}{dx^{-\alpha}}f(x)=\frac{1}{\Gamma(\alpha)}\int_c^x f(t)(x-t)^{\alpha-1}dx+\psi_c(x),\end{equation}
	where $\psi_c(x)$ was an arbitrary ``complimentary'' function meant to eliminate the ambiguity in the lower integration limit \cite{pimentel2017fractional}.
	
	Let it first be recognized that the distributional differintegral does exactly this: eliminates the ambiguity of the lower integration limit. Secondly, however, the distributional differintegral actually proposes a method to solve for Riemann's $\psi_c(x)$ complimentary function.
	
	An interpretation of this phenomena may be seen using inverse derivatives. Since
	\begin{equation}\int_{x_0}^x f(t)dt=\frac{d^{-1}}{dx^{-1}}f(x)-\frac{d^{-1}}{dx^{-1}}f(x)\Big|_{x=x_0}=\text{\LARGE\S}_x^1 f(x)-\bigg(\text{\LARGE\S}_x^1 f(x)\bigg)\bigg|_{x=x_0}.\end{equation}
	
	Repetition of this process, and Cauchy's formula for repeated integration states
	\begin{equation}\frac{1}{(n-1)!}\int_{x_0}^x f(t) (x-t)^{n-1} dt=\frac{d^{-n}}{dx^{-n}}f(x)-\sum_{k=0}^{n-1} c_k x^{n-k-1}=\text{\LARGE\S}_x^n f(x)-\sum_{k=0}^{n-1} c_k x^{n-k-1},\end{equation}
	where $c_k$ is given recursively by the formula
	
	\begin{equation}c_{k}=\frac{d^{-k}}{dx^{-k}}f(x)\Big|_{x=x_0}-\sum_{j=0}^{k-1}c_j x_0^{k-j-1}=\bigg(\text{\LARGE\S}_x^k f(x)\bigg)\bigg|_{x=x_0}-\sum_{j=0}^{k-1}c_j x_0^{k-j-1},\end{equation}
	so that each step reflects the prior evaluated at $x_0$.
	
	Taking the distributional differintegral of both sides to reveal $f(x)$ one recovers the Riemann-Liouville definition
	\begin{equation}\lim_{m\to0^+}\frac{1}{\Gamma(m)}\int_{x_0}^x f(t) (x-t)^{m-1} dt=f(x)-\sum_{k=0}^{n-1} c_k \emptyset^{(k)}(x),\end{equation}
	or rather
	
	\begin{equation}f(x)=\lim_{m\to0^+}\frac{1}{\Gamma(m)}\int_{x_0}^x f(t) (x-t)^{m-1} dt+\sum_{k=0}^{n-1} c_k \emptyset^{(k)}(x).\end{equation}
	
	Since $n$ was arbitrarily large, and since $\emptyset^{(n)}(x)=0$ \textit{almost everywhere} for $n\in\mathbb{N}\cup\{0\}$, one may take the limit as $n\to\infty$, recovering the formulas
	
	\begin{equation}\lim_{m\to0^+}\frac{1}{\Gamma(m)}\int_{x_0}^x f(t) (x-t)^{m-1} dt=f(x)-\sum_{k=0}^{\infty} c_k \emptyset^{(k)}(x),\end{equation}
	and
	
	\begin{equation}f(x)=\lim_{m\to0^+}\frac{1}{\Gamma(m)}\int_{x_0}^x f(t) (x-t)^{m-1} dt+\sum_{k=0}^{\infty} c_k \emptyset^{(k)}(x).\end{equation}
	
	This formula gives the necessary ``constants of integration'' for repeated integration of integer order, and similarly keeps the Zero Functions equal to zero \textit{almost everywhere} for derivatives of integer order.
	
	The issues arise when taking derivatives of non-integer order; that is, fractional derivatives. Observe now that
	
	\begin{equation}f^{(\alpha)}(x)=\frac{1}{\Gamma(-\alpha)}\int_{x_0}^x f(t) (x-t)^{-\alpha-1} dt+\sum_{k=0}^{\infty}c_k \emptyset^{(k+\alpha)}(x),\end{equation}
	or
	
	\begin{equation}\frac{1}{\Gamma(-\alpha)}\int_{x_0}^x f(t) (x-t)^{-\alpha-1} dt=f^{(\alpha)}(x)-\sum_{k=0}^{\infty}c_k \emptyset^{(k+\alpha)}(x).\end{equation}
	
	It will be shown later that fractional (non-integer) derivatives of Zero Functions are \textit{not} equal to zero \textit{almost everywhere} but rather have \textit{nonzero} values everywhere. Because of this, integral definitions of fractional derivatives will \textit{almost always} contain polynomials of infinite degree. This also accounts for why many of these integrals do not converge.
	
	Herein is proposed a solution to Riemann's $\psi_c(x)$ complimentary function. That is
	
	\begin{equation}\psi_c(x)=\sum_{k=0}^{\infty}c_k \emptyset^{(k+\alpha)}(x).\end{equation}
	
	Sadly, this complimentary function was eliminated from the Riemann-Liouville definition due to Laurent's work in 1884.
	
	\section{Specific distributional differintegrals}
	
	Herein are provided distributional differintegrals of common functions of a single real variable. Notice that the operator takes a real-valued function and creates a complex-valued function (though still of a single real variable). Computation for these cases of functions of a single real variable mimic that of the Riemann-Liouville derivative. Therefore they have been omitted for brevity.
	
	Before beginning, it is important to note that the distributional differintegral is linear as shown above. Thus, one must only apply it to portions of functions separated by addition/subtraction. Translations of the independent variable are allowed in all distributional differintegrals. That is, if $\widehat{T}_{z_0}$ is an independent variable translation operator such that $\widehat{T}_{z_0} f(z)=f(z-z_0)$, then
	
	\begin{equation}\text{\LARGE\S}_z ^\alpha \widehat{T}_{z_0}[f(z)]=\widehat{T}_{z_0}\Big[\text{\LARGE\S}_z ^\alpha f(z) \Big].\end{equation}
	This is a result true for any convolution with a distribution, as seen in \cite{sacks2017techniques} exercise 6.14.
	
	In the following results, $\Omega$ is chosen such that $\text{\Large\S}_z^\alpha f(z)\big|_{\inf\Omega}=0$, as this applies when solving most differential equations and allows for better convergence of test functions.
	
	\subsection{Monomials}
	
	The distributional differintegral of a generic monomial function is
	
	\begin{equation}\text{\LARGE\S}_z ^\alpha z^n=\frac{\Gamma(1+n)}{\Gamma(1+n+\alpha)}z^{n+\alpha},\end{equation}
	where $n\in\mathbb{C}\backslash\mathbb{Z}^-$ (or in the special case of Zero Functions, $n\in\mathbb{Z}^-$ and the function is scaled by $\frac{1}{\Gamma(1+n)}$, forcing the scalars to cancel when the distributional differintegral is applied).
	
	Now observe the existence of the Zero Function. Of course expanded from our earlier definition, the Zero Function may now be defined as
	
	\begin{equation}\emptyset(z)=\text{\LARGE\S}_z ^{-1}z^0=\frac{1}{\Gamma(0)}z^{-1}.\end{equation}
	
	Likewise, the $\alpha^{th}$ derivative of the Zero Function may be defined as
	
	\begin{equation}\emptyset^{(\alpha)}(z)=\text{\LARGE\S}_z ^{-\alpha} \emptyset(z)=\text{\LARGE\S}_z ^{-\alpha} \frac{1}{\Gamma(0)}z^{-1}=\frac{1}{\Gamma(0)}\cdot\frac{\Gamma(0)}{\Gamma(-\alpha)}z^{-1-\alpha}=\frac{1}{\Gamma(-\alpha)}z^{-1-\alpha}.\end{equation}
	
	This is fascinating as it implies that non-integer derivatives of Zero Functions are not zero, but notice just as well that every integer derivative of Zero Functions is indeed zero.
	
	It is well-known that the inverse Laplace Transform behaves as
	
	\begin{equation}\mathscr{L}^{-1}[s^n,t]=\frac{1}{2\pi i}\lim_{T\to\infty}\int_{\gamma-iT}^{\gamma+iT}s^n e^{st}ds=\frac{1}{\Gamma(-n)}t^{-1-n}.\end{equation}
	Considering the Zero Function was defined $\emptyset^{(\alpha)}(z)=\frac{1}{\Gamma(-\alpha)}z^{-1-\alpha}$, the result should look quite familiar. Indeed this implies
	
	\begin{equation}\mathscr{L}[\emptyset^{(\alpha)}(t),s]=s^\alpha.\end{equation}
	
	Notice as well, that one could redefine the definition of the distributional differintegral utilizing the Zero Function. Observe,
	
	\begin{equation}\emptyset^{(\alpha)}=\frac{1}{\Gamma(-\alpha)}z^{-1-\alpha} \hspace{0.25 in}\text{or}\hspace{0.25 in} \emptyset^{(-\alpha)}=\frac{1}{\Gamma(\alpha)}z^{\alpha-1},\end{equation}
	implying
	
	\begin{equation}\text{\LARGE\S}_z ^\alpha f(z)=\frac{1}{\Gamma(\alpha)}\bigg(\frac{d^{-1}}{d\zeta^{-1}} f(\zeta)(z-\zeta)^{\alpha-1}\Big\vert_{\zeta=z}\bigg)=\bigg(\frac{d^{-1}}{d\zeta^{-1}} f(\zeta)\emptyset^{(-\alpha)}(z-\zeta)\Big\vert_{\zeta=z}\bigg),\end{equation}
	or
	
	\begin{equation}\text{\LARGE\S}_z ^{-\alpha} f(z)=\bigg(\frac{d^{-1}}{d\zeta^{-1}} f(\zeta)\emptyset^{(\alpha)}(z-\zeta)\Big\vert_{\zeta=z}\bigg).\end{equation}
	This certainly should not be a surprise as the Zero Function is literally \textit{defined} as a result of the distributional differintegral of $z^0$.
	
	\subsection{Dirac delta function}
	The distributional differintegral of the Dirac delta function is
	
	\begin{equation}\text{\LARGE\S}_z ^\alpha \delta(z)=H(z)\frac{z^{\alpha-1}}{\Gamma(\alpha)}=H(z)\emptyset^{(-\alpha)}(z).\end{equation}
	
	This, of course, extends nicely to the distributional differintegral of the Heaviside step function,
	
	\begin{equation}\text{\LARGE\S}_z ^\alpha H(z)=H(z)\frac{z^{\alpha}}{\Gamma(\alpha+1)}=H(z)\emptyset^{(-\alpha-1)}(z),\end{equation}
	
	The answers come as no surprise, since it is well-known that the antiderivatives of the Heaviside step function are the antiderivatives of the monomials multiplied by the Heaviside step function. The first derivative of the Heaviside step function is the Dirac delta function, so it would follow that the Zero Function multiplied by the Heaviside step function is also the Dirac delta function.
	
	One may also rewrite the distributional differintegral as a convolution with a fractional Dirac delta function. Notice that
	
	\begin{equation}\text{\LARGE\S}_z^{-\alpha} f(z)=\int_{\mathbb{R}}f(\zeta)\emptyset^{(\alpha)}(z-\zeta)H(z-\zeta)d\zeta=\int_{\mathbb{R}}f(\zeta)\delta^{(\alpha)}(z-\zeta)d\zeta,\end{equation}
	where the $\alpha^{th}$ fractional derivative of the Dirac Delta Function is denoted as $\delta^{(\alpha)}(z)=\text{\LARGE\S}_z^{-\alpha}\delta(z)$.
	
	\textbf{Remark.} This is where the title of the manuscript arises. The operator is simply a fractionalization of the distributional derivatives and inverse derivatives of the Dirac delta function. Along with the sifting property and integration by parts, it results that \textit{a fractional differintegral may be based on the convolution of the fractionalized distributional derivative of the Dirac delta function}.
	
	\subsection{Exponential Function}
	
	The distributional differintegral of a generic exponential function is
	
	\begin{equation}\text{\LARGE\S}_z ^\alpha e^{\lambda z}=\lambda^{-\alpha}e^{\lambda z}.\end{equation}	
	This, coupled with the linearity of this operator, implies the distributional differintegral of the $\sin$, $\cos$, $\sinh$, and $\cosh$ functions.
	
	\begin{equation}\text{\LARGE\S}_z ^\alpha \sin(\lambda z)=|\lambda|^{-\alpha} \sin\Big(\lambda z -\alpha\frac{\pi}{2}\Big).\end{equation}
	
	\begin{equation}\text{\LARGE\S}_z ^\alpha \cos(\lambda z)=|\lambda|^{-\alpha} \cos\Big(\lambda z -\alpha\frac{\pi}{2}\Big).\end{equation}
	
	\begin{equation}\text{\LARGE\S}_z ^\alpha \sinh(\lambda z)=\bigg(\frac{\lambda^{-\alpha} e^{\lambda z}-(-\lambda)^{-\alpha}e^{-\lambda z}}{2}\bigg).\end{equation}
	
	\begin{equation}\text{\LARGE\S}_z ^\alpha \cosh(\lambda z)=\bigg(\frac{\lambda^{-\alpha} e^{\lambda z}+(-\lambda)^{-\alpha}e^{-\lambda z}}{2}\bigg).\end{equation}
	
	\textbf{Remark.} The Gamma function's definition follows from a specific case of the exponential function's distributional differintegral.
	
	Observe the following:
	
	\begin{equation}\text{\LARGE\S}_z^\alpha e^z=\frac{1}{\Gamma(\alpha)}\int_{-\infty}^{z}e^\zeta(z-\zeta)^{\alpha-1}H(z-\zeta)d\zeta=e^z.\end{equation}
	Thus,
	
	\begin{equation}\Big[\text{\LARGE\S}_z^\alpha e^z\Big] \Big|_{z=0}=\Big[\frac{1}{\Gamma(\alpha)}\int_{-\infty}^{z}e^\zeta(z-\zeta)^{\alpha-1}H(z-\zeta)d\zeta\Big]\Big|_{z=0}=[e^z] \big|_{z=0}.\end{equation}
	This implies,
	
	\begin{equation}\frac{1}{\Gamma(\alpha)}\int_{-\infty}^{0}e^\zeta(0-\zeta)^{\alpha-1}H(0-\zeta)d\zeta=e^0=1.\end{equation}
	
	One may rearrange this as
	\begin{equation}\Gamma(\alpha)=\int_{-\infty}^{0}e^\zeta (0-\zeta)^{\alpha-1}H(0-\zeta)d\zeta=\int_{0}^{\infty}\zeta^{\alpha-1} e^\zeta d\zeta,\end{equation}
	which may be recognized as the most-common integral definition of the Gamma function.
	
	\subsection{Natural Logarithm}
	
	The distributional differintegral of a the natural logarithm function is
	
	\begin{equation}\text{\LARGE\S}_z ^\alpha \ln(\lambda z)=z^\alpha \frac{\ln(z)+\ln(\lambda)-\gamma-\psi(1+\alpha)}{\Gamma(1+\alpha)},\end{equation}
	where $\gamma\approx0.57721...$, the Euler-Mascheroni constant, and $\psi(1+\alpha)=\frac{\Gamma'(1+\alpha)}{\Gamma(1+\alpha)}$, the digamma function.
	
	Note here, that $\ln(\lambda z)=\ln(z)+\ln(\lambda)=\ln(z)+\ln(\lambda)z^0$, and indeed the linearity exists,
	
	\begin{equation}\text{\LARGE\S}_z ^\alpha \Big[\ln(z)+\ln(\lambda)\Big]=z^\alpha \frac{\ln(z)-\gamma-\psi(1+\alpha)}{\Gamma(1+\alpha)}+\ln(\lambda)\cdot \frac{1}{\Gamma(1+\alpha)}z^\alpha=\text{\LARGE\S}_z ^\alpha \ln(\lambda z).\end{equation}
	
	It is in this distributional differintegral that the strongest advantage of the distributional differintegral is found. When attempting to create a distributional differintegral, there was prior doubt as to how functions such as $z^{-n}$ would be treated. On one hand, they were monomials, but on the other, they were the derivatives of logarithms. There was fear that a ``perfect'' definition would never be created because of this discrepancy.
	
	What may be seen, however, is that the ``monomial'' distributional differintegral of $z^{-n}$ is found in the form $\frac{\Gamma(1-n)}{\Gamma(1-n-\alpha)}z^{-n-\alpha}$, where at integer-valued derivatives, the result is one of the derivatives of the Zero Function. Likewise, the distributional differintegral of the ``logarithmic'' version of $z^{-n}$ is found in the form $z^\alpha \frac{\ln(z)+\ln(\lambda)-\gamma-\psi(1+\alpha)}{\Gamma(1+\alpha)}$, where at integer-valued derivatives, the result keeps the functional form $z^{-n}$ and is not killed by an infinite denominator (as a result of the digamma function). Interestingly, the non-integer derivatives of this class of functions include the natural logarithm again.
	
	\subsection{The Polylogarithm, and thus the Riemann Zeta Function, are specific cases of the distributional differintegral}
	
	Observe the following distributional differintegral:
	
	\begin{equation}\text{\LARGE\S}_z ^\alpha \frac{1}{e^{-z}-1}={Li}_\alpha(e^z),\end{equation}
	where $Li_\alpha (e^z)$ is the polylogarithm function of base $\alpha$.
	
	Testing the distributional differintegral on this function was no coincidence. Note the similarity of the Riemann Zeta Function's integral definition,
	
	\begin{equation}\zeta(s)=\frac{1}{\Gamma(s)}\int_0^\infty \frac{x^{s-1}}{e^x-1}dx,\end{equation}
	and the integral definition of the distributional differintegral. A quick rearrangement of the Zeta Function shows,
	
	\begin{equation}
	\begin{split}
	\zeta(s)&=\frac{-1}{\Gamma(s)}\int_\infty^0 \frac{1}{e^x -1}(x-0)^{s-1}H(x-0)dx=\frac{1}{\Gamma(s)}\int_{-\infty}^{0} \frac{1}{e^{-x}-1}(0-x)^{s-1}H(0-x)dx\\
	&=\frac{1}{\Gamma(s)}\bigg(\int_{-\infty}^{t} \frac{1}{e^{-x}-1}(t-x)^{s-1}H(t-x)dx\bigg)\bigg|_{t=0} =\bigg(\text{\LARGE\S}_t ^s \frac{1}{e^{-t}-1}\bigg)\bigg|_{t=0}.
	\end{split}
	\end{equation}
	
	Indeed it may be seen from above that
	
	\begin{equation}\bigg(\text{\LARGE\S}_z ^\alpha \frac{1}{e^{-z}-1}\bigg)\bigg|_{z=0}={Li}_\alpha(e^z)\big|_{z=0}={Li}_\alpha(1)=\zeta(\alpha).\end{equation}
	
	This also supports the fact that the polylogarithm may be defined as the repeated integral of itself. Here the distributional differintegral power $\alpha$ as the base of the polylogarithm implies even further that the definition holds true.
	
	\subsection{Product of Monomial and Exponential}
	
	The distributional differintegral of a monomial-exponential product is
	
	\begin{equation}\text{\LARGE\S}_z ^\alpha z^n e^{\lambda z}=\alpha\frac{\Gamma(1+n)}{\Gamma(1+\alpha)}(\lambda z)^{n+\alpha}\cdot {}_1 F_1 (1+n;1+\alpha;\lambda z)-\frac{\alpha\Gamma(-\alpha)}{\Gamma(1-n-\alpha)}\lambda^{-n-\alpha}\cdot {}_1 F_1 (1-\alpha;1-n-\alpha;\lambda z),\end{equation}
	where ${}_1 F_1 (1+n;1+\alpha;\lambda z)$ and ${}_1 F_1 (1-\alpha;1-n-\alpha;\lambda z)$ are Kummer confluent hypergeometric functions.
	
	While the distributional differintegral becomes complicated very quickly for products of functions (no doubt this arises from a generalization of product rules and integration by parts), the product of the exponential and monomial is a very commonly-used product, notably in the fractional Laplace transform.
	
	\section{Distributional Differintegral Transforms}
	
	\subsection{Fractional Laplace Transform}
	
	It is possible to generalize the Laplace Transform to fractional values using the distributional differintegral. This is done as follows, utilizing the operator $\mathscr{L}^{(\alpha)}[f(t),s]$ as the $\alpha^{th}$ power of the Laplace Transform.
	
	For $0\leq\alpha\leq1$,
	
	\begin{equation}\mathscr{L}^{(\alpha)}[f(t),s]=e^{-i\pi\alpha} \bigg(e^{t s}\text{\LARGE\S}_t^\alpha f(t)e^{-s t}\bigg)\bigg|_{t=(1-\alpha)s}.\end{equation}
	
	One may quickly see that
	
	\begin{equation}\mathscr{L}^{(0)}[f(t),s]=f(s),\end{equation}
	and with slightly more effort
	
	\begin{equation}\mathscr{L}^{(1)}[f(t),s]=\mathscr{L}[f(t),s]=F(s),\end{equation}
	where $F(s)$ is the Laplace Transform of $f(t)$.
	
	\textbf{Proof.}
	
	\begin{equation}
	\begin{split}
	\mathscr{L}^{(0)}[f(t),s]&=e^{-i\pi\cdot 0} \bigg(e^{t s}\text{\LARGE\S}_t^0 f(t)e^{-s t}\bigg)\bigg|_{t=s}\\
	&=e^{s^2}\text{\LARGE\S}_s^0 f(s) e^{-s^2}=e^{s^2-s^2}f(s)\\
	&=f(s),
	\end{split}
	\end{equation}
	and
	
	\begin{equation}
	\begin{split}
	\mathscr{L}^{(1)}[f(t),s]&=e^{-i\pi\cdot 1} \bigg(e^{t s}\text{\LARGE\S}_t^1 f(t)e^{-s t}\bigg)\bigg|_{t=0}=-\bigg(e^{t s}\int_{\mathbb{R}} f(\tau)e^{-s \tau}(t-\tau)^0 H(t-\tau)d\tau\bigg)\bigg|_{t=0}\\
	&=-e^{0}\int_{\mathbb{R}} f(\tau)e^{-s \tau}H(-\tau)d\tau=\int_\mathbb{R} f(\tau)e^{-s\tau}H(\tau)d\tau=\int_0^\infty f(\tau)e^{-s\tau}d\tau\\
	&=\mathscr{L}[f(t),s]=F(s).
	\end{split}
	\end{equation}
	
	Notice that one may use this definition without including the evaluation of $t=(1-\alpha)s$. However, when using this evaluation, one arrives at a function of a single independent variable (instead of \textit{two} independent variables) for fractional values of $\alpha$. There is little intuition as to what the interpretation of a result in this manner would imply phsyically.
	
	\subsection{Differintegral Fourier Transform}
	
	As the fractional Fourier transform already exists, a different name must be used for the transform resulting from the distributional differintegral applied to the Fourier transform. Since the Fourier Transform is a specific case of the bilateral Laplace Transform, it is possible to construct a Fourier Transform from a special sum of the Laplace Transform above. This is done as follows, utilizing the operator $\mathscr{F}^{(\alpha)}[f(t),\omega]$ as the $\alpha^{th}$ power of the Fourier Transform.
	
	For $0\leq\alpha\leq1$,
	
	\begin{equation}
	\begin{split}
	\mathscr{F}^{(\alpha)}[f(t),\omega]&=\frac{1}{4}\bigg(\frac{2}{\pi}\bigg)^\frac{\alpha}{2} \bigg(e^{i\omega t}\text{\LARGE\S}_t^\alpha f(t)e^{-i\omega t}\bigg)\bigg|_{t=(1-\alpha)\omega}+\frac{1}{4}\bigg(\frac{2}{\pi}\bigg)^\frac{\alpha}{2}\bigg(e^{-i\omega t}\text{\LARGE\S}_t^\alpha f(-t)e^{i\omega t}\bigg)\bigg|_{t=(1-\alpha)\omega}\\
	&+\frac{1}{4}\bigg(\frac{2}{\pi}\bigg)^\frac{\alpha}{2}\bigg(e^{i\omega t}\text{\LARGE\S}_t^\alpha f(t)e^{-i\omega t}\bigg)\bigg|_{t=(\alpha-1)\omega}+\frac{1}{4}\bigg(\frac{2}{\pi}\bigg)^\frac{\alpha}{2}\bigg(e^{-i\omega t}\text{\LARGE\S}_t^\alpha f(-t)e^{i\omega t}\bigg)\bigg|_{t=(\alpha-1)\omega}.
	\end{split}
	\end{equation}
	
	Because of the relationship of the Laplace Transform to the Fourier Transform, the results for $\alpha=0$ and $\alpha=1$ are equivalent in their parts to the above fractional Laplace Transform at $\alpha=0$ and $\alpha=1$.
	
	\begin{equation}\mathscr{F}^{(0)}[f(t),\omega]=f(\omega),\end{equation}
	and
	
	\begin{equation}\mathscr{F}^{(1)}[f(t),\omega]=\widehat{f}(\omega),\end{equation}
	where $\widehat{f}(\omega)$ is the Fourier Transform of $f(t)$.
	
	\section{Fractional Distributional Differintegral Equations}
	
	It is helpful to recognize the result that
	
	\begin{equation}\mathscr{L}\Big[\text{\LARGE\S}_t ^\alpha f(t),s\Big]=s^{-\alpha}F(s)=s^{-\alpha}\mathscr{L}[f(t),s],\end{equation}
	or in the form of derivatives, with notation $f^{(\alpha)}(t)=\text{\Large\S}_t^{-\alpha}f(t)$,
	
	\begin{equation}\mathscr{L}\Big[\text{\LARGE\S}_t ^{-\alpha} f(t),s\Big]=\mathscr{L}[f^{(\alpha)}(t),s]=s^{\alpha}F(s)=s^{\alpha}\mathscr{L}[f(t),s].\end{equation}
	
	In the past, initial conditions would be introduced in the $s$-space as derivatives of the function evaluated at $0$. In the case of fractional differential equations, these initial conditions are introduced by means of Zero Functions. Indeed the inclusion of Zero Functions reconciles a very important question regarding fractional differential equations. One knows that for an $n^{th}$ order differential equation, one needs $n$ boundary conditions to solve it. How many initial conditions are necessary for an $n=\frac{1}{2}$ order differential equation? Using a Laplace Transform (with the ``sum of derivatives'' subtracted from the derivative) one ends up with notation such as
	
	\begin{equation}\sum_{k=1}^{1/2} s^{n-k}f^{(k-1)}(0),\end{equation}
	which is quite ambiguous. Indeed, summation notation is built with integer values in mind. Past fractional differential equations have been forced to use $\lceil\alpha\rceil$ boundary conditions. This works quite well, but sometimes requires more boundary conditions than necessary.
	
	\subsection{Generalization of the nonlinear Volterra integral equation}
	
	Because the Volterra integral equation allows one to rewrite a differential equation as an integral equation, one may add initial conditions (and later boundary conditions) to a fractional differential equation.
	
	Suppose a fractional differential equation is of the form
	
	\begin{equation}y^{(\alpha)}(x)=f\big(x,y(x)\big).\end{equation}
	Then past research \cite{diethelm2002analysis} insists its equivalent Volterra integral equation is
	
	\begin{equation}y(x)=\sum_{k=0}^{\lceil\alpha\rceil-1} y^{(k)}(0)\frac{x^k}{k!}+\frac{1}{\Gamma(\alpha)}\int_0^x f\big(\tau,y(\tau)\big)(x-\tau)^{\alpha-1}d\tau.\end{equation}
	
	Interestingly, Zero Functions and the distributional differintegral allow even further generalization of this formula. That is, for the sequence of $n$ arbitrary initial conditions, $y^{(\alpha_k)}(0)$ (note these may be fractional derivatives of $f$), and $k=1,2,...,n$.
	
	\begin{equation}y(x)=\sum_{k=1}^{n}y^{(\alpha_k)}(0)\frac{x^{\alpha_k}}{\Gamma(1+\alpha_k)}+\text{\LARGE\S}_x^\alpha f\big(x,y(x)\big).\end{equation}
	
	Taking the distributional differintegral (to the $-\alpha^{th}$ power) of both sides, one arrives at a new version of the fractional differential equation
	
	\begin{equation}y^{(\alpha)}(x)=\sum_{k=1}^n y^{(\alpha_k)}(0)\frac{x^{\alpha_k-\alpha}}{\Gamma(1+\alpha_k-\alpha)}+f\big(x,y(x)\big),\end{equation}
	or with Zero Functions,
	
	\begin{equation}y^{(\alpha)}(x)=\sum_{k=1}^n y^{(\alpha_k)}(0)\emptyset^{(\alpha-\alpha_k-1)}(x)+f\big(x,y(x)\big).\end{equation}
	
	Just so, if all of the values of $(\alpha-\alpha_k)\in\mathbb{N}$, as one sees in the current theory of differential equations, then this equation is equivalent \textit{almost everywhere} to the original.
	
	It is even possible to extend this theory to arbitrary boundary conditions and not initial conditions, though it makes for a more difficult equation to solve. Instead of the initial conditions in the Volterra integral equations, one may leave the constant as an arbitrary $c_k$, but then there must be enough constants left at the end to satisfy all initial conditions. This leaves the following result, for the sequence of $n$ arbitrary \textit{boundary} conditions, $y^{(\alpha_k)}(x_k)$ (note these may also be fractional derivatives of $f$), and $k=1,2,...,n$:
	
	\begin{equation}y(x)=\sum_{k=1}^{n}c_k x^{\alpha-k}+\text{\LARGE\S}_x^\alpha f\big(x,y(x)\big),\end{equation}
	or in differential equation form with Zero Functions,
	
	\begin{equation}y^{(\alpha)}(x)=\sum_{k=1}^n c_k \emptyset^{(k-1)}(x)+f\big(x,y(x)\big).\end{equation}
	
	In this case, the differential equation is solved without computing the values of the constants first, then the values of each constant is solved by means of a system of equations.
	
	\subsection{Contraction Mapping Theorem in Banach Space}
	
	Since the space of functions $X^{-c}(\Omega)$ was assumed to have a subspace dense in norm $\|\cdot\|_{X^{-c}(\Omega)}$, the space of must be at the very least a normed linear space. If the space happens to be complete (as it is in most cases), then it is a Banach space.
	
	Thus suppose $X^{-c}(\Omega)$ is complete. Then for a bounded linear operator $T$, if $T^n$ is a contraction ($\|T^n\|<1$)for some power $n$, there exists a fixed point $f\in X^{-c}(\Omega)$ such that $T(f)=f$ \cite{sacks2017techniques}. This theorem becomes useful in solving differential and integral equations as finding the fixed point is equivalent to solving the equation.
	
	In $X^{-c}(\Omega)$, with an operator defined as $F(x)=Tx+b$ for some $b\in X^{-c}(\Omega)$, if $\|T\|<1$ then the fixed point solution to the equation is
	\begin{equation}x=\sum_{j=0}^\infty T^j b.\end{equation}
	
	Observe that for all $j\in\mathbb{C}$ (and thus all $j\in\mathbb{N}$)
	
	\begin{equation}\bigg(\text{\LARGE\S}_z^\alpha\bigg)^j=\text{\LARGE\S}_z^{j\alpha}\end{equation}
	as seen from the existence of the index law. This implies that for any equation of the form
	
	\begin{equation}u(z)=F\big(u(z)\big)=\text{\LARGE\S}_z^\alpha u(z)+f(z)\end{equation}
	where $f\in X^{-c}(\Omega)$ and $\text{\Large\S}_z^\alpha$ a contraction, one may solve the equation to arbitrary approximation with
	
	\begin{equation}u(z)=\sum_{j=0}^\infty \text{\LARGE\S}_z^{j\alpha}f(z).\end{equation}
	When $\text{\Large\S}_z^\alpha f(z)=g(\alpha,z)$ this allows for even easier computation with
	
	\begin{equation}u(z)=\sum_{j=0}^\infty g(j\alpha,z).\end{equation}
	
	\section{Conclusion}
	
	In conclusion, a consistent definition for the distributional differintegral was established. This definition allows extension of differentiation and inverse differentiation to all complex powers. The definition only behaves, however, if one eliminates the so-called ``constant of integration'' by means of using Zero Functions. These Zero Functions have interesting properties and indeed seem to hold the secret to the distributional differintegral.
	
	It is most common to compute the distributional differintegral by performing a definite integral with an introduced Heaviside step function, but can also be formed using the ``inverse derivatives'' briefly introduced. These equivalent definitions satisfy the four properties necessary for a fractional derivative given by \cite{ortigueira2015fractional}.
	
	Compared with past distributional fractional derivatives, these results may  be a reconciliation of what was proposed in \cite{stojanovic2011generalized}, which stated that the Riemann-Liouville derivative operating in a distributional sense does not produce an integer-valued distributional derivative. Including the notion of Zero Functions in the results from \cite{stojanovic2011generalized} gives precisely the integer-valued distributional derivatives.
	
	The specific distributional differintegrals of functions agree with the most-commonly used fractional derivatives of many functions, with the distributional differintegral of monomials appearing as Riemann-Liouville's definition, and the distributional differintegral of an exponential function appearing as that of the Grunwald-Letnikov definition. The cases of the natural logarithm, polylogarithm, and others appear slightly different from many previous definitions. A special case of the polylogarithm shows the emergence of the Riemann-Zeta function, while a special case of the exponential function shows the emergence of the gamma function.
	
	There is also a fractional Laplace Transform that may be introduced using the fractional integral portion of the distributional differintegral. This may be generalized to a differintegral Fourier Transform.
	
	A small number of fractional differential equations may be solved with the Laplace Transform of the distributional differintegral. Inclusion of the Zero Functions in these equations rectifies a slight issue in determining how many initial conditions are necessary to solve fractional differential equations. Other fractional differential equations may be solved using a generalized Volterra integral equation. One may also simplify the fixed-point solution to a contraction mapping theorem in the Banach space of analytic functions using the distributional differintegral.
	
	While the properties of the distributional differintegral in regards to differential equations were only briefly examined, the future is very bright for all of the areas covered in this paper. It is quite likely that the surface has only begun to be scratched on the power and potential of the distributional differintegral.
	
	\section*{Acknowledgments}
	
	The author would like to thank Dr. Darin J. Ulness and Dr. Douglas R. Anderson for their immense help in revising the manuscript as well as for stimulating conversation about the topic.\\
	
	\textbf{Declaration of Interest.} The author declares no financial and personal relationships with other people or organizations that could inappropriately influence this work.

\end{document}